\documentclass{amsart}%
\usepackage{amsmath}
\usepackage{graphicx}
\usepackage{amsfonts}
\usepackage[all]{xy}
\usepackage{amscd}
\usepackage{amssymb}%
\newtheorem{theorem}{Theorem}
\theoremstyle{plain}

\newtheorem{corollary}{Corollary}

\newtheorem{definition}{Definition}

\newtheorem{lemma}{Lemma}

\newtheorem{openproblem}{Open problem}
\newtheorem{proposition}{Proposition}
\newtheorem{remark}{Remark}

\numberwithin{equation}{section}
\def\rotl#1{\rotdimen=\ht#1\advance\rotdimen by\dp#1   \hbox to\rotdimen{\vbox to\wd#1{\vskip\wd#1\rotstart{270 rotate}   \box#1\vss}\hss}\rotfinish}
\begin{document}
\title[Collective properties of \ concave Young-functions]{Bijections and metric spaces induced by some collective properties of concave Young-functions}
\author{N. K. Agbeko}
\address{Institute of Mathematics\\
University of Miskolc\\
H-3515 Miskolc--Egyetemv\'aros\\
Hungary}
\email{matagbek@uni-miskolc.hu}
\date{March 24th, 2006}
\subjclass{Primary 26A06, 54E35, 26A42; Secondary 11J83, 28A25, 47H10}
\keywords{concave Young-functions, dense sets, metric spaces, bijections}

\begin{abstract}
For each ${\small b\in\left(  0,\,\infty\right)  }$ we intend to generate a
decreasing sequence of subsets $\left(  \mathcal{Y}_{b}^{\left(  n\right)
}\right)  \subset Y_{\mathrm{conc}}$ depending on $b$\ such that whenever
$n\in\mathbb{N}$, then $\mathcal{A}\cap\mathcal{Y}_{b}^{\left(  n\right)  }%
$\ is dense in $\mathcal{Y}_{b}^{\left(  n\right)  }$ and the following four
sets $\mathcal{Y}_{b}^{\left(  n\right)  }$, $\mathcal{Y}_{b}^{\left(
n\right)  }\backslash\left(  \mathcal{A}\cap\mathcal{Y}_{b}^{\left(  n\right)
}\right)  $, $\mathcal{A}\cap\mathcal{Y}_{b}^{\left(  n\right)  }$
and\ $\mathcal{Y}_{\mathrm{conc}}$\ are pairwise equinumerous. Among others we
also show that if $f$ is any measurable function on a measure space $\left(
\Omega,\mathcal{F},\lambda\right)  $ and $p\in\left[  1,\infty\right)  $ is an
arbitrary number then the quantities $\left\Vert f\right\Vert _{L^{p}}$ and
$\sup_{\Phi\in\widetilde{\mathcal{Y}_{\mathrm{conc}}}}\left(  \Phi\left(
1\right)  \right)  ^{-1}\left\Vert \Phi\circ\left\vert f\right\vert
\right\Vert _{L^{p}}$ are equivalent, in the sense that they are both either
finite or infinite at the same time.

\end{abstract}
\maketitle

\section{Introduction}

We know that concave functions play major roles in many branches of
mathematics for instance probability theory (\cite{BURK1973}, \cite{GARS1973},
\cite{MOGY1981}, say), interpolation theory (cf. \cite{TRIEB1978}, say),
weighted norm inequalities (cf. \cite{GARFRAN1985}, say), and functions spaces
(cf. \cite{SINN2002}, say), as well as in many other branches of sciences. In
the line of \cite{BURK1973}, \cite{GARS1973} and \cite{MOGY1981}, the present
author also obtained in martingale theory some results in connection with
certain collective properties or behaviors of concave Young-functions (cf.
\cite{AGB1986}, \cite{AGB1989}). The study presented in \cite{AGB2005} was
mainly motivated by the question why strictly concave functions possess so
many properties, worth to be characterized using appropriate tools that await
to be discovered.

We say that a function $\Phi:\left[  0,\,\infty\right)  \rightarrow\left[
0,\,\infty\right)  $ belongs to the set $\mathcal{Y}_{\mathrm{conc}}$ (and is
referred to as a concave Young-function) if and only if it admits the integral
representation%
\begin{equation}
\Phi\left(  x\right)  =\int\nolimits_{0}^{x}\varphi\left(  t\right)  dt,
\label{id1}%
\end{equation}
(where $\varphi:\left(  0,\,\infty\right)  \rightarrow\left(  0,\,\infty
\right)  $ is a right-continuous and decreasing function such that it is
integrable on every finite interval $\left(  0,\,x\right)  $) and $\Phi\left(
\infty\right)  =\infty$. It is worth to note that every function in
$\mathcal{Y}_{\mathrm{conc}}$ is strictly concave.

We will remind some results obtained so far in \cite{AGB2005}.

We shall say that a concave Young-function $\Phi$ satisfies the
\emph{density-level property} if $A_{\Phi}\left(  \infty\right)  <\infty$,
where $A_{\Phi}\left(  \infty\right)  :=\int\nolimits_{1}^{\infty}%
\frac{\varphi\left(  t\right)  }{t}dt$. All the concave Young-functions
possessing the density-level property will be grouped in a set $\mathcal{A}$.

In Theorems \ref{theo1} and \ref{theo2} (cf. \cite{AGB2005}), we showed that
the composition of any two concave Young-functions satisfies the density-level
property if and only if at least one of them satisfies it. These two theorems
show that concave Young-functions with the density-level property behave like
left and right ideal with respect to the composition operation.

We also proved (\cite{AGB2005}, Lemma 5, page 12) that if\ $\Phi\in
\mathcal{Y}_{\mathrm{conc}}$, then\ there are constants $C_{\Phi}>0$\ and
$B_{\Phi}\geq0$ such that%
\[
A_{\Phi}\left(  \infty\right)  -B_{\Phi}\leq%
{\displaystyle\int_{0}^{\infty}}
\frac{\Phi\left(  t\right)  }{\left(  t+1\right)  ^{2}}dt\leq C_{\Phi}%
+A_{\Phi}\left(  \infty\right)  .
\]
This led us to the idea to search for a Lebesgue measure (described here
below) with respect to which every concave Young-function turns out to be
square integrable (\cite{AGB2005}, Lemma 6, page 13), i.e. $\mathcal{Y}%
_{\mathrm{conc}}\subset L^{2}:=L^{2}\left(  \left[  0,~\infty\right)
,~\mathcal{M},~\mu\right)  $, where $\mathcal{M}$ is a $\sigma$-algebra (of
$\left[  0,~\infty\right)  $) containing the Borel sets and $\mu
:\mathcal{M}\rightarrow\left[  0,~\infty\right)  $ is a Lebesgue measure
defined by $\mu\left(  \left[  0,~x\right)  \right)  =\frac{1}{3}\left(
1-\frac{1}{\left(  x+1\right)  ^{3}}\right)  $ for all $x\in\left[
0,~\infty\right)  $. The mapping $d:L^{2}\times L^{2}\rightarrow\left[
0,~\infty\right)  $, defined by%
\begin{equation}
\operatorname*{d}\left(  f,g\right)  =\sqrt{%
{\displaystyle\int_{\left[  0,~\infty\right)  }}
\left(  f-g\right)  ^{2}d\mu}=\sqrt{%
{\displaystyle\int_{0}^{\infty}}
\frac{\left(  f\left(  x\right)  -g\left(  x\right)  \right)  ^{2}}{\left(
x+1\right)  ^{4}}dx}, \label{dist}%
\end{equation}
is known to be a semi-metric.

Further on, we proved in (\cite{AGB2005}, Theorem 8, page 16) that
$\mathcal{A}$\ is a dense set in $\mathcal{Y}_{\mathrm{conc}}$.

Throughout this communication $\Phi_{\operatorname{id}}$ will denote the
identity function defined on the half line $\left[  0,\text{ }\infty\right)  $
and we write $\left\Vert \Phi\right\Vert :=\sqrt{%
{\displaystyle\int_{\left[  0,\text{ }\infty\right)  }}
\Phi^{2}d\mu}$\ whenever $\Phi\in\mathcal{Y}_{\mathrm{conc}}$.

We intend to generate a decreasing sequence of subsets $\left(  \mathcal{Y}%
_{b}^{\left(  n\right)  }\right)  \subset\mathcal{Y}_{\mathrm{conc}}$
depending on $b$\ such that whenever $n\in\mathbb{N}$, then $\mathcal{A}%
\cap\mathcal{Y}_{b}^{\left(  n\right)  }$\ is dense in $\mathcal{Y}%
_{b}^{\left(  n\right)  }$ and the following four sets $\mathcal{Y}%
_{b}^{\left(  n\right)  }$, $\mathcal{Y}_{b}^{\left(  n\right)  }%
\backslash\left(  \mathcal{A}\cap\mathcal{Y}_{b}^{\left(  n\right)  }\right)
$, $\mathcal{A}\cap\mathcal{Y}_{b}^{\left(  n\right)  }$ and\ $\mathcal{Y}%
_{\mathrm{conc}}$\ are pairwise equinumerous. We shall also prove that the two
pairs $\left(  \mathcal{Z}^{\ast\left(  n\right)  },\operatorname*{dist}%
\right)  $\ and $\left(  \mathcal{Z}^{\left(  n\right)  },\operatorname*{dist}%
\right)  $ are metric spaces, where $\mathcal{Z}^{\ast\left(  n\right)
}=\left\{  \mathcal{Y}_{b}^{\left(  n\right)  }:{\small b\in\left(
0,\,\infty\right)  }\right\}  $ and $\mathcal{Z}^{\left(  n\right)  }=\left\{
\mathcal{A}_{b}^{\left(  n\right)  }:{\small b\in\left(  0,\,\infty\right)
}\right\}  $ for each $n\in\mathbb{N}$ and the distance between any two sets
$\mathcal{F}$ and $\mathcal{G}$ in $\mathcal{Y}_{\mathrm{conc}}$ being defined
by
\begin{align*}
\operatorname*{dist}\left(  \mathcal{F},\mathcal{G}\right)   &  :=\sup\left\{
\inf\left\{  \operatorname{d}\left(  \Phi,\Psi\right)  :\Psi\in\mathcal{G}%
\right\}  :\Phi\in\mathcal{F}\right\} \\
&  =\sup\left\{  \inf\left\{  \operatorname{d}\left(  \Phi,\Psi\right)
:\Phi\in\mathcal{F}\right\}  :\Psi\in\mathcal{G}\right\}  .
\end{align*}
We show in the last section that if $f$ is any measurable function on a
measure space $\left(  \Omega,\mathcal{F},\lambda\right)  $ and $p\in\left[
1,\infty\right)  $ is an arbitrary number then the quantities $\left\Vert
f\right\Vert _{L^{p}}$ and $\sup_{\Phi\in\widetilde{\mathcal{Y}_{\mathrm{conc}%
}}}\left(  \Phi\left(  1\right)  \right)  ^{-1}\left\Vert \Phi\circ\left\vert
f\right\vert \right\Vert _{L^{p}}$ are equivalent, in the sense that they are
both either finite or infinite at the same time, where $\widetilde
{\mathcal{Y}_{\mathrm{conc}}}$ is a proper subset of $\mathcal{Y}%
_{\mathrm{conc}}$.We then use this subset to express the value of $\left\Vert
f\right\Vert _{L^{p}}$ whenever $\left\Vert f\right\Vert _{L^{p}}<\infty$.

\section{Bijections between subsets of $\mathcal{Y}_{\mathrm{conc}}$}

We first anticipate that there are as many elements in each of the sets
$\mathcal{A}$ and $\mathcal{Y}_{\mathrm{conc}}\backslash\mathcal{A}$\ as there
exist in $\mathcal{Y}_{\mathrm{conc}}$, showing how broad the set of concave
Young-functions possessing the density-level property and its complement
really are.

\begin{theorem}
\label{theo1}The sets $\mathcal{A}$, $\mathcal{Y}_{\mathrm{conc}}$ and
$\mathcal{Y}_{\mathrm{conc}}\backslash\mathcal{A}$ are pairwise equinumerous.
\end{theorem}

\begin{proof}
We first show that there is a bijection between $\mathcal{A}$ and
$\mathcal{Y}_{\mathrm{conc}}$. In fact, since $\mathcal{A}$ is a proper subset
of $\mathcal{Y}_{\mathrm{conc}}$ there is an injection from $\mathcal{A}$\ to
$\mathcal{Y}_{\mathrm{conc}}$, as a matter of fact, the identity mapping from
$\mathcal{A}$ into $\mathcal{Y}_{\mathrm{conc}}$ will do. Fix any number
$\alpha\in\left(  0,\,1\right)  $ and define the mapping $S_{\alpha
}:\mathcal{Y}_{\mathrm{conc}}\rightarrow\mathcal{A}$ by $S_{\alpha}\left(
\Phi\right)  =\Phi^{\alpha}$. We point out that this mapping exists in virtue
of Theorem 2 in \cite{AGB2005}. It is not hard to see that $S_{\alpha}$ is an
injection. Then the Schr\"{o}der-Bernstein theorem entails that there exists a
bijection between $\mathcal{A}$ and $\mathcal{Y}_{\mathrm{conc}}$. To complete
the proof it is enough to show that there is a bijection between $\mathcal{A}$
and $\mathcal{Y}_{\mathrm{conc}}\backslash\mathcal{A}$. In fact, fix
arbitrarily some $\Phi\in\mathcal{Y}_{\mathrm{conc}}\backslash\mathcal{A}$ and
define the function $h_{\Phi}:\mathcal{A}\rightarrow\mathcal{Y}_{\mathrm{conc}%
}\backslash\mathcal{A}$ by $h_{\Phi}\left(  \Delta\right)  =\Delta+\Phi$.
Obviously, $h_{\Phi}$ is an injection. Now, fix any $\Delta\in\mathcal{A}$ and
define the function $f_{\Delta}:\mathcal{Y}_{\mathrm{conc}}\backslash
\mathcal{A}\rightarrow\mathcal{A}$ by $f_{\Delta}\left(  \Phi\right)
=\Delta\circ\Phi$. We point out that this function always exists due to
Theorem \ref{theo2} in \cite{AGB2005}. It is not difficult to show that
$f_{\Delta}$ is an injection if we take into account that $\Delta$ is an
invertible function. Consequently, the Schr\"{o}der-Bernstein theorem
guarantees the existence of a bijection between $\mathcal{A}$ and
$\mathcal{Y}_{\mathrm{conc}}\backslash\mathcal{A}$. Therefore, we can conclude
on the validity of the argument.
\end{proof}

Write $\mathcal{A}_{b}:=\left\{  \Phi\in\mathcal{A}:\Phi\left(  b\right)
=b\right\}  $ and $\mathcal{Y}_{b}:=\left\{  \Phi\in\mathcal{Y}_{\mathrm{conc}%
}:\Phi\left(  b\right)  =b\right\}  $ for every number $b\in\left(
0,\,\infty\right)  $.

Let us denote by $\mathcal{Z}:=\left\{  \mathcal{A}_{b}:b\in\left(
0,\,\infty\right)  \right\}  $ and $\mathcal{Z}^{\ast}:=\left\{
\mathcal{Y}_{b}:b\in\left(  0,\,\infty\right)  \right\}  $.

It is obvious that $\mathcal{A}_{b}\subset\mathcal{Y}_{b}$ for every number
$b\in\left(  0,\,\infty\right)  $ and $\mathcal{Z}\cap\mathcal{Z}^{\ast
}=\varnothing$.

\begin{lemma}
\label{lem1}For every number $b\in\left(  0,\,\infty\right)  $ the identities
$\mathcal{A}_{b}=\left\{  \frac{b\Phi}{\Phi\left(  b\right)  }:\Phi
\in\mathcal{A}\right\}  $ and $\mathcal{Y}_{b}=\left\{  \frac{b\Phi}%
{\Phi\left(  b\right)  }:\Phi\in\mathcal{Y}_{\mathrm{conc}}\right\}  $ hold true.
\end{lemma}

\begin{proof}
Pick any function $\Psi\in\mathcal{A}_{b}$. Then $\Psi\in\mathcal{A}$ and
$\Psi\left(  b\right)  =b$, so that $\Psi=\frac{b\Psi}{\Psi\left(  b\right)
}\in\left\{  \frac{b\Phi}{\Phi\left(  b\right)  }:\Phi\in\mathcal{A}\right\}
$, i.e. $\mathcal{A}_{b}\subset\left\{  \frac{b\Phi}{\Phi\left(  b\right)
}:\Phi\in\mathcal{A}\right\}  $. To show the reverse inclusion consider any
function $\Psi\in\left\{  \frac{b\Phi}{\Phi\left(  b\right)  }:\Phi
\in\mathcal{A}\right\}  $. Then necessarily there must exist some $\Phi
\in\mathcal{A}$ such that $\Psi=\frac{b\Phi}{\Phi\left(  b\right)  }$. It is
obvious that $\Psi\in\mathcal{A}$ and $\Psi\left(  b\right)  =b$, i.e.
$\Psi\in\mathcal{A}_{b}$. Hence, $\left\{  \frac{b\Phi}{\Phi\left(  b\right)
}:\Phi\in\mathcal{A}\right\}  \subset\mathcal{A}_{b}$. These two inclusions
yield that $\mathcal{A}_{b}=\left\{  \frac{b\Phi}{\Phi\left(  b\right)  }%
:\Phi\in\mathcal{A}\right\}  $. The proof of identity $\mathcal{Y}%
_{b}=\left\{  \frac{b\Phi}{\Phi\left(  b\right)  }:\Phi\in\mathcal{Y}%
_{\mathrm{conc}}\right\}  $\ can be similarly carried out.
\end{proof}

\begin{definition}
\label{def1}A proper subset $\mathcal{G}$ of $\mathcal{A}$ is said to be
maximally bounded if each of the sets $\mathcal{G}$ and $\mathcal{A}%
\backslash\mathcal{G}$ is equinumerous with $\mathcal{A}$, i.e. there is a
bijection between $\mathcal{A}$\ and $\mathcal{G}$, and $\operatorname*{diam}%
(\mathcal{G})<\infty$, where $\operatorname*{diam}(\mathcal{G}):=\sup\left\{
\operatorname*{d}\left(  \Phi_{1},\Phi_{2}\right)  :\Phi_{1},~\Phi_{2}%
\in\mathcal{G}\right\}  $ is the diameter of $\mathcal{G}$.
\end{definition}

We note that Definition \ref{def1} makes sense for the two reasons here below.

On the one hand we assert that $\operatorname*{diam}(\mathcal{A})=\sup\left\{
\operatorname*{d}\left(  \Phi_{1},\Phi_{2}\right)  :\Phi_{1},~\Phi_{2}%
\in\mathcal{A}\right\}  =\infty$. In fact, fix some $\Phi\in\mathcal{A}$
and\ define a sequence $\left(  \Phi_{n}\right)  \subset\mathcal{Y}%
_{\mathrm{conc}}$ by $\Phi_{2n}=4n\Phi$ and $\Phi_{2n-1}=\left(  2n-1\right)
\Phi$, $n\in\mathbb{N}$. It is clear that $\left(  \Phi_{n}\right)
\subset\mathcal{A}$ and $\operatorname*{d}\left(  \Phi_{2n},\Phi
_{2n-1}\right)  =\left(  2n+1\right)  \left\Vert \Phi\right\Vert $,
$n\in\mathbb{N}$. Hence, $\operatorname*{diam}(\mathcal{A})=\infty$.

On the other hand the set $\left\{  \left(  \Phi\left(  1\right)  \right)
^{-1}\Phi:\Phi\in\mathcal{Y}_{\mathrm{conc}}\right\}  $ is of finite diameter.
In fact for any $\Phi$, $\Psi\in\mathcal{Y}_{\mathrm{conc}}$ we have, via
Lemma \ref{lem3} in \cite{AGB2005}, that
\[
\operatorname*{d}\left(  \left(  \Phi\left(  1\right)  \right)  ^{-1}%
\Phi,\left(  \Psi\left(  1\right)  \right)  ^{-1}\Psi\right)  \leq\left\Vert
\left(  \Phi\left(  1\right)  \right)  ^{-1}\Phi\right\Vert +\left\Vert
\left(  \Psi\left(  1\right)  \right)  ^{-1}\Psi\right\Vert \leq2\left\Vert
S\right\Vert <\infty.
\]

Let us define two relations $\mathrm{\bot}\subset\mathcal{A}\times\mathcal{A}%
$\ and $\mathrm{\bot}^{\ast}\subset\mathcal{Y}_{\mathrm{conc}}\times
\mathcal{Y}_{\mathrm{conc}}$ as follows:

\begin{enumerate}
\item We say that $\left(  \Phi,\Psi\right)  \in\mathrm{\bot}$, where $\left(
\Phi,\Psi\right)  \in\mathcal{A}\times\mathcal{A}$, (and write $\Phi
\mathrm{\bot}\Psi$) if and only if there is some constant $c\in\left(
0,\,\infty\right)  $\ such that $\Psi\left(  x\right)  =c\Phi\left(  x\right)
$\ for all $x\in\left(  0,\,\infty\right)  $.

\item We say that $\left(  \Phi,\Psi\right)  \in\mathrm{\bot}^{\ast}$, where
$\left(  \Phi,\Psi\right)  \in\mathcal{Y}_{\mathrm{conc}}\times\mathcal{Y}%
_{\mathrm{conc}}$, (and write $\Phi\mathrm{\bot}^{\ast}\Psi$) if and only if
there is some constant $c\in\left(  0,\,\infty\right)  $\ such that
$\Psi\left(  x\right)  =c\Phi\left(  x\right)  $\ for all $x\in\left(
0,\,\infty\right)  $.
\end{enumerate}

It is not hard to see that $\mathrm{\bot}$ and $\mathrm{\bot}^{\ast}$ are
equivalence relations on $\mathcal{A}$ and $\mathcal{Y}_{\mathrm{conc}}$
respectively, i.e. they are reflexive, symmetric and transitive. Their
corresponding equivalence classes are respectively%
\begin{align*}
p_{\mathrm{\bot}}\left(  \Psi\right)   &  :=\left\{  \Phi:\Phi\in
\mathcal{A}\text{ and }\Phi\mathrm{\bot}\Psi\right\}  ,\text{ }\Psi
\in\mathcal{A}\\
p_{\mathrm{\bot}^{\ast}}\left(  \Delta\right)   &  :=\left\{  \Phi:\Phi
\in\mathcal{Y}_{\mathrm{conc}}\text{ and }\Phi\mathrm{\bot}^{\ast}%
\Delta\right\}  ,\text{ }\Delta\in\mathcal{Y}_{\mathrm{conc}}%
\end{align*}
and their respective induced factor (or quotient) sets can be given by%
\begin{align*}
\mathcal{A}/\mathrm{\bot}  &  :=\left\{  \mathcal{C}:\mathcal{C}%
\subset\mathcal{A}\text{ and }\mathcal{C}=p_{\mathrm{\bot}}\left(
\Psi\right)  \text{ for some }\Psi\in\mathcal{A}\right\}  ,\\
\mathcal{Y}_{\mathrm{conc}}/\mathrm{\bot}^{\ast}  &  :=\left\{  \mathcal{C}%
:\mathcal{C}\subset\mathcal{Y}_{\mathrm{conc}}\text{ and }\mathcal{C}%
=p_{\mathrm{\bot}^{\ast}}\left(  \Delta\right)  \text{ for some }\Delta
\in\mathcal{Y}_{\mathrm{conc}}\right\}
\end{align*}

One can easily verify that for all $\Psi\in\mathcal{A}$ and $\Delta
\in\mathcal{Y}_{\mathrm{conc}}$\ the equivalence classes\ $p_{\mathrm{\bot}%
}\left(  \Psi\right)  $\ and $p_{\mathrm{\bot}^{\ast}}\left(  \Delta\right)
$\ are of continuum size or magnitude.

\begin{theorem}
\label{theo2}Let $b\in\left(  0,\,\infty\right)  $ be any fixed
number.\newline\textbf{Part I.} Define the mapping $f:\mathcal{A}%
\rightarrow\mathcal{A}_{b}$ by $f\left(  \Phi\right)  =\frac{b}{\Phi\left(
b\right)  }\Phi$. Then there is a unique mapping\ $g:\mathcal{A}/\mathrm{\bot
}\rightarrow\mathcal{A}_{b}$ for which the diagram%
\begin{equation}
\xymatrix{\mathcal{A} \ar[r]^{p_{\mathrm{\bot}}} \ar[dr]_f & \mathcal{A}/\mathrm{\bot} \ar[d]^g\\ ~ & \mathbf{b} }
\label{com3}%
\end{equation}
commutes \emph{(}i.e. $f=g\circ p_{\mathrm{\bot}}$\emph{)} and moreover, the
mapping $g$ is a bijection.\newline\textbf{Part II.} Define the mapping
$f^{\ast}:\mathcal{Y}_{\mathrm{conc}}\rightarrow\mathcal{Y}_{b}$ by $f^{\ast
}\left(  \Delta\right)  =\frac{b}{\Delta\left(  b\right)  }\Delta$. Then there
is a unique mapping\ $g^{\ast}:\mathcal{Y}_{\mathrm{conc}}/\mathrm{\bot}%
^{\ast}\rightarrow\mathcal{Y}_{b}$ for which the diagram%
\begin{equation}
\xymatrix{\mathcal{Y}_{\mathrm{conc}} \ar[r]^{p_{\mathrm{\bot}^{\ast}}} \ar[dr]_{f^{\ast}} & \mathcal{Y}_{\mathrm{conc}}/\mathrm{\bot}^{\ast} \ar[d]^{g^{\ast}}\\ ~ & \mathbf{b}^{\ast} }
\label{com4}%
\end{equation}
commutes \emph{(}i.e. $f^{\ast}=g^{\ast}\circ p_{\mathrm{\bot}^{\ast}}%
$\emph{)} and moreover, the mapping $g^{\ast}$ is a bijection.
\end{theorem}

We point out that the proof of Theorem \ref{theo2} is obvious.

\begin{proposition}
\label{prop1}Let $b\in\left(  0,\,\infty\right)  $ be an arbitrarily fixed
number.\newline\textbf{Part I.} There is a bijection between $\mathcal{Y}_{b}$
and $\mathcal{Y}_{\mathrm{conc}}$.\newline\textbf{Part II.} There is a
bijection between $\mathcal{A}_{b}$ and $\mathcal{A}$.
\end{proposition}

\begin{proof}
We shall only show the first part because the other case can be similarly
proved. To this end, write $\mathcal{Y}_{bb}:=\left\{  b\Phi:\Phi
\in\mathcal{Y}_{\mathrm{conc}}\right\}  $. We note that $\mathcal{Y}_{bb}$ and
$\mathcal{Y}_{\mathrm{conc}}$ are equinumerous for the reasons that
$\mathcal{Y}_{bb}\subset\mathcal{Y}_{\mathrm{conc}}$ and the function
$F:\mathcal{Y}_{\mathrm{conc}}\rightarrow\mathcal{Y}_{bb}$, defined by
$F\left(  \Phi\right)  =b\Phi$, can be easily shown to be an injection. Thus
it will be enough to prove that $\mathcal{Y}_{bb}$ and $\mathcal{Y}_{b}$ are
equinumerous. In fact, consider the function $Q:\mathcal{Y}_{b}\rightarrow
\mathcal{Y}_{bb}$ defined by $Q\left(  \frac{b}{\Phi\left(  b\right)  }%
\Phi\right)  =b\Phi$. We shall just point out that function $Q$ can be easily
shown to be a bijection, which ends the proof.
\end{proof}

\begin{corollary}
\label{cor1}Let $b\in\left(  0,\,\infty\right)  $ be arbitrary. Then the
following six sets $\mathcal{A}$, $\mathcal{A}_{b}$, $\mathcal{Y}_{b}$,
$\mathcal{A}/\mathrm{\bot}$, $\mathcal{Y}_{\mathrm{conc}}/\mathrm{\bot}^{\ast
}$and $\mathcal{Y}_{\mathrm{conc}}$ are equinumerous.
\end{corollary}

\begin{proof}
We note that $\mathcal{A}$ and $\mathcal{Y}_{\mathrm{conc}}$ are equinumerous
(by Theorem \ref{theo1}) and, by Theorem 2, $\mathcal{A}/\mathrm{\bot}$ and
$\mathcal{A}_{b}$ are equinumerous. On the other hand $\mathcal{A}$ and
$\mathcal{A}_{b}$ are equinumerous as well as $\mathcal{Y}_{b}$ and
$\mathcal{Y}_{\mathrm{conc}}$\ are (by Proposition \ref{prop1}). Thus
$\mathcal{A}_{b}$ and $\mathcal{Y}_{\mathrm{conc}}$ are equinumerous.
Therefore, as $\mathcal{Y}_{b}$ and $\mathcal{Y}_{\mathrm{conc}}/\mathrm{\bot
}^{\ast}$ are equinumerous\ (by Theorem 2), we can conclude on the validity of
the argument.
\end{proof}

\begin{remark}
\label{rem1}Let $b_{1}$ and $b_{2}\in\left(  0,\,\infty\right)  $ be two
arbitrary distinct numbers. Then $\mathcal{A}_{b_{1}}\cap\mathcal{A}_{b_{2}}$
and $\mathcal{Y}_{b_{1}}\cap\mathcal{Y}_{b_{2}}$ are empty sets.
\end{remark}

\begin{remark}
\label{rem2}Let $b_{1}$ and $b_{2}\in\left(  0,\,\infty\right)  $ be two
arbitrary distinct numbers. Then $\mathcal{A}_{b_{1}}\cup\mathcal{A}_{b_{2}%
}\notin\mathcal{Z}$ and $\mathcal{Y}_{b_{1}}\cup\mathcal{Y}_{b_{2}}%
\notin\mathcal{Z}^{\ast}$.
\end{remark}

\begin{remark}
\label{rem3}Fix arbitrarily a number $b\in\left(  0,\,\infty\right)  $. Then
it is easily seen that the function $h_{b}:\left[  0,~\infty\right)
\rightarrow\left[  0,~\infty\right)  $, defined by $h_{b}\left(  x\right)
=x+b$, is square integrable with respect to measure $\mu$ and, moreover,
$C_{b}:=%
{\displaystyle\int_{0}^{\infty}}
\frac{\left(  h_{b}\left(  x\right)  \right)  ^{2}}{\left(  x+1\right)  ^{4}%
}dx=\frac{1}{3}\left(  b^{2}+b+1\right)  <\infty$.
\end{remark}

\begin{remark}
\label{rem4}If $\Phi\in\mathcal{Y}_{b}$, then $\Phi\left(  x\right)  \leq
h_{b}\left(  x\right)  $ for all $x\in\left[  0,\,\infty\right)  $.
\end{remark}

\begin{proof}
Fix any $\Phi\in\mathcal{Y}_{b}$. As $\Phi$ is a concave function its graph
must lie below the tangent of equation $y=\varphi\left(  b\right)  \left(
x-b\right)  +b$ at point $\left(  b,b\right)  $ since $\Phi\left(  b\right)
=b$. Consequently, for all $x\in\left[  0,\,\infty\right)  $ we have:%
\begin{align*}
\Phi\left(  x\right)   &  \leq\varphi\left(  b\right)  \left(  x-b\right)
+b\leq\varphi\left(  b\right)  x+b=b\varphi\left(  b\right)  \frac{x}{b}+b\\
&  \leq\Phi\left(  b\right)  \frac{x}{b}+b=h_{b}\left(  x\right)  .
\end{align*}

\end{proof}

\begin{proposition}
\label{prop2}Let $b\in\left(  0,\,\infty\right)  $ be any number. Then
$\mathcal{Y}_{b}$ is of finite diameter.
\end{proposition}

\begin{proof}
Let $b\in\left(  0,\,\infty\right)  $ be the source of $\mathcal{Y}_{b}%
\in\mathcal{Z}^{\ast}$. We need to prove that $\mathcal{Y}_{b}$ has a finite
diameter. In fact, consider two arbitrary functions $\Phi_{1}$, $\Phi_{2}%
\in\mathcal{Y}_{b}$. Then%
\[
\operatorname{d}\left(  \Phi_{1},\Phi_{2}\right)  =\left\Vert \Phi_{1}%
-\Phi_{2}\right\Vert \leq\left\Vert \Phi_{1}\right\Vert +\left\Vert \Phi
_{2}\right\Vert \leq\sqrt{2C_{b}}\text{,}%
\]
via Remarks \ref{rem4} and \ref{rem3}. Therefore,%
\[
\operatorname*{diam}\left(  \mathcal{Y}_{b}\right)  :=\sup\left\{
\operatorname{d}\left(  \Phi_{1},\Phi_{2}\right)  :\Phi_{1},~\Phi_{2}%
\in\mathcal{Y}_{b}\right\}  \leq\sqrt{2C_{b}}<\infty.
\]

\end{proof}

\begin{theorem}
\label{theo3}Let $b\in\left(  0,\,\infty\right)  $ be any number. Then
$\mathcal{Y}_{b}$ is maximally bounded.
\end{theorem}

\begin{proof}
We just point out that the proof follows from the conjunction of both
Propositions \ref{prop2} and \ref{prop1}.
\end{proof}

In the sequel $H_{\left[  0,~1\right]  }$ will stand for the collection of all
finite sequences $\left(  t_{1},~\ldots,~t_{k}\right)  \subset\left[
0,~1\right]  $ such that $t_{1}+\ldots+t_{k}=1$.

For any fixed $b\in\left(  0,\,\infty\right)  $ and every counting number
$n\in\mathbb{N}$ write $\operatorname*{X}\limits_{i=1}^{n}\mathcal{A}_{b}$
(resp. $\operatorname*{X}\limits_{i=1}^{n}\mathcal{Y}_{b}$) for the $n$-fold
Descartes product of $\mathcal{A}_{b}$ (resp. $\mathcal{Y}_{b}$).

\begin{center}
For $n=1$ let us set $\mathcal{A}_{b}^{\left(  1\right)  }=\mathcal{A}_{b}$,
$\mathcal{Y}_{b}^{\left(  1\right)  }=\mathcal{Y}_{b}$ and whenever $n\geq2$,
write $\mathcal{Y}_{b}^{CO\left(  n\right)  }=\left\{  \Delta_{1}\circ
\Delta_{2}\circ\ldots\circ\Delta_{n}:\left(  \Delta_{1},\Delta_{2}%
,\ldots,\Delta_{n}\right)  \in\operatorname*{X}\limits_{i=1}^{n}%
\mathcal{Y}_{b}\right\}  $,
\[
\mathcal{A}_{b}^{CO\left(  n\right)  }=\left\{  \Phi_{1}\circ\ldots\circ
\Phi_{n}:\left(  \Phi_{1},\ldots,\Phi_{n}\right)  \in\operatorname*{X}%
\limits_{i=1}^{n}\mathcal{Y}_{b}\text{ and }\Phi_{j}\in\mathcal{A}_{b}\text{
for some index }j\right\}  ,
\]
$\mathcal{Y}_{b}^{\left(  n\right)  }=\left\{  \sum_{i=1}^{k}t_{i}\Delta
_{i}:\Delta_{1},~\Delta_{2},~\ldots,~\Delta_{k}\in\mathcal{Y}_{b}^{CO\left(
n\right)  },~\left(  t_{1},~\ldots~t_{k}\right)  \in H_{\left[  0,~1\right]
}\right\}  $, $\mathcal{A}_{b}^{\left(  n\right)  }=\left\{  \sum_{i=1}%
^{k}t_{i}\Phi_{i}:\Phi_{1},~\Phi_{2},~\ldots,~\Phi_{k}\in\mathcal{A}%
_{b}^{CO\left(  n\right)  },~\left(  t_{1},~\ldots~t_{k}\right)  \in
H_{\left[  0,~1\right]  }\right\}  $.
\end{center}

\noindent Further, for $n=1$ write $\mathcal{Z}^{\left(  1\right)
}=\mathcal{Z}$, $\mathcal{Z}^{\ast\left(  1\right)  }=\mathcal{Z}^{\ast}$ and,
for $n\in\mathbb{N}\backslash\left\{  1\right\}  $ write $\mathcal{Z}^{\left(
n\right)  }:=\left\{  \mathcal{A}_{b}^{\left(  n\right)  }:b\in\left(
0,\,\infty\right)  \right\}  $ and $\mathcal{Z}^{\ast\left(  n\right)
}:=\left\{  \mathcal{Y}_{b}^{\left(  n\right)  }:b\in\left(  0,\,\infty
\right)  \right\}  $.

\begin{remark}
\label{rem5}For any pair of numbers $n\in\mathbb{N}$ and $b\in\left(
0,\,\infty\right)  $ the set $\mathcal{A}_{b}^{\left(  n\right)  }$ is a
proper subset of $\mathcal{Y}_{b}^{\left(  n\right)  }$.
\end{remark}

\begin{remark}
\label{rem6}For any pair of numbers $n\in\mathbb{N}$ and $b\in\left(
0,\,\infty\right)  $ we have $\mathcal{A}_{b}^{\left(  n\right)  }%
\subset\mathcal{A}_{b}^{\left(  1\right)  }=\mathcal{A}_{b}$.
\end{remark}

We point out that Remark \ref{rem6} is a direct consequent of Theorem 2 in
\cite{AGB2005}, page 6.

\begin{remark}
\label{rem7}Let $b\in\left(  0,\,\infty\right)  $, $n\in\mathbb{N}$ and $k\geq
n$ be arbitrary numbers. Then \newline\emph{(1) }$\Phi_{1}\circ\Phi_{2}%
\circ\ldots\circ\Phi_{k}\in\mathcal{A}_{b}^{CO\left(  n\right)  }$ whenever
$\Phi_{1},~\Phi_{2},~\ldots,~\Phi_{k}\in\mathcal{Y}_{b}^{\left(  1\right)  }$
and $\Phi_{j}\in\mathcal{A}_{b}^{\left(  1\right)  }$ for some index
$j\in\left\{  1,~\ldots,~k\right\}  $\newline\emph{(2) }$\Delta_{1}\circ
\Delta_{2}\circ\ldots\circ\Delta_{k}\in\mathcal{Y}_{b}^{CO\left(  n\right)  }$
whenever $\Delta_{1},~\Delta_{2},~\ldots,~\Delta_{k}\in\mathcal{Y}%
_{b}^{\left(  1\right)  }$.
\end{remark}

\begin{proof}
Note that $\Phi_{1}\circ\Phi_{2}\circ\ldots\circ\Phi_{k}=\Phi_{1}\circ\Phi
_{2}\circ\ldots\circ\Phi_{n-1}\circ\Psi_{1}$ and $\Delta_{1}\circ\Delta
_{2}\circ\ldots\circ\Delta_{n-1}\circ\Psi_{2}$, where $\Psi_{1}=\Phi_{n}%
\circ\Phi_{n+1}\circ\ldots\circ\Phi_{k}$ and $\Psi_{2}=\Delta_{n}\circ
\Delta_{n+1}\circ\ldots\circ\Delta_{k}$. From this simple observation the
result easily follows.
\end{proof}

From Remark \ref{rem7} the following result can be easily derived, since it
implies that $\mathcal{A}_{b}^{CO\left(  n+1\right)  }$ is a proper subset of
$\mathcal{A}_{b}^{CO\left(  n\right)  }$ and, $\mathcal{Y}_{b}^{CO\left(
n+1\right)  }$ is also a proper subset of $\mathcal{Y}_{b}^{CO\left(
n\right)  }$.

\begin{lemma}
\label{lem2}Let $b\in\left(  0,\,\infty\right)  $ and $n\in\mathbb{N}$ be
arbitrary numbers. Then the following two assertions are valid.\newline%
\emph{(1)} The set $\mathcal{A}_{b}^{\left(  n+1\right)  }$ is a proper subset
of $\mathcal{A}_{b}^{\left(  n\right)  }$.\newline\emph{(2) }The set
$\mathcal{Y}_{b}^{\left(  n+1\right)  }$ is a proper subset of $\mathcal{Y}%
_{b}^{\left(  n\right)  }$.
\end{lemma}

\begin{theorem}
\label{theo4}For any fixed pair of numbers $n\in\mathbb{N}$ and $b\in\left(
0,\,\infty\right)  $, the two sets $\mathcal{A}_{b}^{\left(  n\right)  }$ and
$\mathcal{A}_{b}$ are equinumerous.
\end{theorem}

\begin{proof}
Throughout the proof we shall fix any counting number $n\in\mathbb{N}$. We
first note that the identity function $I_{\operatorname{id}}:\mathcal{A}%
_{b}^{\left(  n\right)  }\rightarrow\mathcal{A}_{b}$ is an injection, since
$\mathcal{A}_{b}^{\left(  n\right)  }\subset\mathcal{A}_{b}$. Next, pick any
$\Delta\in\mathcal{A}_{b}$ and define the function $f_{\Delta}:\mathcal{A}%
_{b}\rightarrow\mathcal{A}_{b}^{\left(  n\right)  }$ by $f_{\Delta}\left(
\Phi\right)  =\underset{\left(  n-1\right)  \text{-fold}}{\underbrace
{\Delta\circ\ldots\circ\Delta}}\circ\Phi$. We show that $f_{\Delta}$ is an
injection. In fact, let $\Phi_{1}$, $\Phi_{2}\in\mathcal{A}_{b}$ be arbitrary
and assume that $f_{\Delta}\left(  \Phi_{1}\right)  =f_{\Delta}\left(
\Phi_{2}\right)  $. Then taking into account that $\Delta$ is an invertible
function we can easily deduce that $\Phi_{1}=\Phi_{2}$, i.e. $f_{\Delta}$ is
an injection. Therefore, the Schr\"{o}der-Bernstein theorem entails that there
is a bijection between $\mathcal{A}_{b}$ and $\mathcal{A}_{b}^{\left(
n\right)  }$. This was to be proved.
\end{proof}

\begin{proposition}
\label{prop3}For any pair of numbers $n\in\mathbb{N}$ and $b\in\left(
0,~\infty\right)  $ the sets $\mathcal{A}_{b}^{\left(  n\right)  }$ and
$\mathcal{Y}_{b}^{\left(  n\right)  }\backslash\mathcal{A}_{b}^{\left(
n\right)  }$ are equinumerous.
\end{proposition}

\begin{proof}
Let $\Phi\in\mathcal{Y}_{b}^{\left(  n\right)  }\backslash\mathcal{A}%
_{b}^{\left(  n\right)  }$ and $\left(  \alpha,\beta\right)  \in H_{\left[
0,~1\right]  }$ be arbitrarily fixed. Define the function $h_{\Phi}^{\left(
\alpha,\beta\right)  }:\mathcal{A}_{b}^{\left(  n\right)  }\rightarrow
\mathcal{Y}_{b}^{\left(  n\right)  }\backslash\mathcal{A}_{b}^{\left(
n\right)  }$ by $h_{\Phi}^{\left(  \alpha,\beta\right)  }\left(
\Delta\right)  =\alpha\Delta+\beta\Phi$. It is clear that $h_{\Phi}^{\left(
\alpha,\beta\right)  }$ is actually an injection. Now, fix any $\Delta
\in\mathcal{A}_{b}^{\left(  n\right)  }$ and define the function $f_{\Delta
}:\mathcal{Y}_{b}^{\left(  n\right)  }\backslash\mathcal{A}_{b}^{\left(
n\right)  }\rightarrow\mathcal{A}_{b}^{\left(  n\right)  }$ by $f_{\Delta
}\left(  \Phi\right)  =\Delta\circ\Phi$. We note that this function always
exists because of the inclusion $\mathcal{A}_{b}^{\left(  n\right)  }%
\subset\mathcal{A}$ and Theorem \ref{theo2} in \cite{AGB2005}. Here too we can
easily check that $f_{\Delta}$ is an injection. Therefore, The
Schr\"{o}der-Bernstein theorem yields the result to be proven.
\end{proof}

\begin{corollary}
For any pair of numbers $n\in\mathbb{N}$ and $b\in\left(  0,~\infty\right)  $
the following five sets $\mathcal{Y}_{b}^{\left(  n\right)  }$, $\mathcal{A}%
_{b}^{\left(  n\right)  }$, $\mathcal{Y}_{b}^{\left(  n\right)  }%
\backslash\mathcal{A}_{b}^{\left(  n\right)  }$, $\mathcal{A}$ and
$\mathcal{Y}_{\mathrm{conc}}$ are pairwise equinumerous.
\end{corollary}

\section{The metrization of sets $\mathcal{Z}^{\left(  n\right)  }$ and
$\mathcal{Z}^{\ast\left(  n\right)  }$}

We shall only deal with the metrization of sets $\mathcal{Z}$ and
$\mathcal{Z}^{\ast}$ since all the results in this section can be easily
extended to the sets $\mathcal{Z}^{\left(  n\right)  }$ and $\mathcal{Z}%
^{\ast\left(  n\right)  }$.

Whenever $\Phi\in\mathcal{Y}_{\mathrm{conc}}$ write $G_{\Phi}:=\left\{
\left(  x,\Phi\left(  x\right)  \right)  :x\in\left(  0,~\infty\right)
\right\}  $ for the graph of $\Phi$ on $\left(  0,~\infty\right)  $ and
$G_{\Phi}^{a||b}:=\left\{  \left(  x,\Phi\left(  x\right)  \right)
:x\in\left[  a,~b\right]  \right\}  $ for the graph of $\Phi$ on the interval
$\left[  a,~b\right]  $ where $a<b$ are any non-negative numbers.

\begin{remark}
\label{rem8}Let $b_{1}$ and $b_{2}\in\left(  0,\,\infty\right)  $ be two
arbitrary distinct numbers. If $b_{1}<b_{2}$, then the following two
assertions hold true:\newline\emph{(1)} For all $\Phi_{1}\in\mathcal{A}%
_{b_{1}}$ and $\Phi_{2}\in\mathcal{A}_{b_{2}}$\ the inequality $\Phi
_{1}\left(  b_{2}\right)  <\Phi_{2}\left(  b_{1}\right)  $\ holds.\newline%
$\emph{(2)}$ For all $\Phi_{1}\in\mathcal{Y}_{b_{1}}$ and $\Phi_{2}%
\in\mathcal{Y}_{b_{2}}$\ the inequality $\Phi_{1}\left(  b_{2}\right)
<\Phi_{2}\left(  b_{1}\right)  $\ holds.
\end{remark}

\begin{proof}
Suppose that $b_{1}<b_{2}$ and fix arbitrarily two functions $\Phi_{1}%
\in\mathcal{Y}_{b_{1}}$ and $\Phi_{2}\in\mathcal{Y}_{b_{2}}$. Obviously,
$\Phi_{1}$ must hit $\Phi_{\operatorname{id}}$\textit{ prior to }$\Phi_{2}$.
Hence, $G_{\Phi_{1}}^{b_{1}||b_{2}}$ lies below $G_{\Phi_{2}}^{b_{1}||b_{2}}$.
But since $G_{\Phi_{1}}^{b_{1}||\infty}$ lies above the graph of the line of
equation $y=b_{1}$ in the interval $\left(  b_{1},\,\infty\right)  $, we have
as an aftermath that $\Phi_{1}\left(  b_{1}\right)  <\Phi_{1}\left(
b_{2}\right)  <\Phi_{2}\left(  b_{1}\right)  $. To end the proof we note that
assertion (2) can be similarly shown.
\end{proof}

The binary relations $\prec$ and $\preceq$\ , defined on $\mathcal{Z}$
respectively by $\mathcal{A}_{b_{1}}\prec\mathcal{A}_{b_{2}}$ if and only if
$\Phi_{1}\left(  b_{2}\right)  <\Phi_{2}\left(  b_{1}\right)  $ for all pairs
$\left(  \Phi_{1},\Phi_{2}\right)  \in\mathcal{A}_{b_{1}}\times\mathcal{A}%
_{b_{2}}$, and by $\mathcal{A}_{b_{1}}\preceq\mathcal{A}_{b_{2}}$ if and only
if $\mathcal{A}_{b_{1}}\prec\mathcal{A}_{b_{2}}$ or $\mathcal{A}_{b_{1}%
}=\mathcal{A}_{b_{2}}$. We point out that The binary relations $\prec$ and
$\preceq$\ can be similarly defined on $\mathcal{Z}^{\ast}$.

We point out that the law of trichotomy is valid on $\left(  \mathcal{Z}%
,\preceq\right)  $ and $\left(  \mathcal{Z}^{\ast},\preceq\right)  $, i.e.
whenever $\left(  \mathcal{A}_{b_{1}},\mathcal{A}_{b_{2}}\right)
\in\mathcal{Z}\times\mathcal{Z}$ or $\left(  \mathcal{A}_{b_{1}}%
,\mathcal{A}_{b_{2}}\right)  \in\mathcal{Z}^{\ast}\times\mathcal{Z}^{\ast}$,
then precisely one of the following holds: $\mathcal{A}_{b_{1}}=\mathcal{A}%
_{b_{2}}$, $\mathcal{A}_{b_{1}}\prec\mathcal{A}_{b_{2}}$, $\mathcal{A}_{b_{2}%
}\prec\mathcal{A}_{b_{1}}$. Hence, we can easily check that $\left(
\mathcal{Z},\preceq\right)  $ and $\left(  \mathcal{Z}^{\ast},\preceq\right)
$ are chains, i.e. they are totally ordered sets.

\begin{theorem}
\label{theo5}The functions $f_{1}:\left(  0,\,\infty\right)  \rightarrow
\mathcal{Z}$ and $f_{2}:\left(  0,\,\infty\right)  \rightarrow\mathcal{Z}%
^{\ast}$, defined respectively by $f_{1}\left(  p\right)  =\mathcal{A}_{p}$
and $f_{2}\left(  p\right)  =\mathcal{Y}_{p}$, are order preserving bijections.
\end{theorem}

\begin{proof}
We show that the function $f_{1}:\left(  0,\,\infty\right)  \rightarrow
\mathcal{Z}$, $f_{1}\left(  p\right)  =\mathcal{A}_{p}$, is an order
preserving bijection. In fact, it is not hard to see via Remark \ref{rem1}
that $f_{1}$ is an injection. Now pick any element $\mathcal{C}\in\mathcal{Z}%
$. Obviously, there must exist some number $p\in\left(  0,\,\infty\right)  $
such that $\mathcal{C}=\mathcal{A}_{p}=f_{1}\left(  p\right)  $, i.e. $f_{1}$
is a surjection. Consequently, $f_{1}$ is a bijection. To end the proof\ of
this part we simply point out that the bijection $f_{1}$ is order preserving
in virtue of Remark \ref{rem5111}. Finally, we note that we can similarly
prove that $f_{2}$ is also an order preserving bijection.
\end{proof}

Since the sets $\left(  \mathcal{Z},\preceq\right)  $ and $\left(
\mathcal{Z}^{\ast},\preceq\right)  $ are chains it is natural to look for a
metric on them. We shall do this in the following two results. But before that
let us recall the definitions of some distances known in the literature (cf.
\cite{KUR1966}, say). If $\Phi\in\mathcal{Y}_{\mathrm{conc}}$\ is any function
and $\mathcal{F}$, $\mathcal{G}\subset\mathcal{Y}_{\mathrm{conc}}$\ are
arbitrary non-empty subsets, then we define the distance from the point $\Phi$
to the set $\mathcal{G}$ by $\rho\left(  \Phi,\mathcal{G}\right)
:=\inf\left\{  \operatorname{d}\left(  \Phi,\Psi\right)  :\Psi\in
\mathcal{G}\right\}  =\inf\left\{  \operatorname{d}\left(  \Psi,\Phi\right)
:\Psi\in\mathcal{G}\right\}  =\rho\left(  \mathcal{G},\Phi\right)  $ and the
distance between the two sets $\mathcal{F}$ and $\mathcal{G}$ by
\begin{align*}
\operatorname*{dist}\left(  \mathcal{F},\mathcal{G}\right)   &  :=\sup\left\{
\inf\left\{  \operatorname{d}\left(  \Phi,\Psi\right)  :\Psi\in\mathcal{G}%
\right\}  :\Phi\in\mathcal{F}\right\} \\
&  =\sup\left\{  \inf\left\{  \operatorname{d}\left(  \Phi,\Psi\right)
:\Phi\in\mathcal{F}\right\}  :\Psi\in\mathcal{G}\right\}  .
\end{align*}

First we find sufficient conditions for which the distance from a point to a
subset (both in $\mathcal{Y}_{\mathrm{conc}}$) should be positive, in order to
guarantee that the distance between two sets in $\mathcal{Y}_{\mathrm{conc}}$
have sense.

\begin{lemma}
\label{lem3}Let $b_{1}$ and $b_{2}\in\left(  0,\,\infty\right)  $ be two
arbitrary distinct numbers. Then $\rho\left(  \mathcal{Y}_{b_{1}},\Phi
_{2}\right)  >0$\ and $\rho\left(  \mathcal{A}_{b_{1}},\Phi_{2}\right)
>0$\ whenever $\Phi_{2}\in\mathcal{Y}_{b_{2}}$.
\end{lemma}

\begin{proof}
It is enough to show that $\rho\left(  \mathcal{A}_{b_{1}},\Phi_{2}\right)
>0$\ whenever $\Phi_{2}\in\mathcal{Y}_{b_{2}}$. In fact, suppose in the
contrary that $\rho\left(  \mathcal{A}_{b_{1}},\Phi_{2}\right)  =0$\ for some
$\Phi_{2}\in\mathcal{Y}_{b_{2}}$. Then there can be extracted some sequence
$\left(  \Delta_{n}\right)  \subset\mathcal{A}_{b_{1}}$ such that
$\operatorname{d}\left(  \Delta_{n},\Phi_{2}\right)  >\operatorname{d}\left(
\Delta_{n+1},\Phi_{2}\right)  $, $n\in\mathbb{N}$, and $\lim_{n\rightarrow
\infty}\operatorname{d}\left(  \Delta_{n},\Phi_{2}\right)  =\rho\left(
\mathcal{A}_{b_{1}},\Phi_{2}\right)  =0$. We point out that this can be done
because of the definition of the infimum. For each $n\in\mathbb{N}$\ let us
set $\Gamma_{n}:=\inf_{k\geq n}\left(  \Delta_{k}-\Phi_{2}\right)  ^{2}%
$.\ Clearly, $\left(  \Gamma_{n}\right)  $ is a non-decreasing sequence of
measurable functions with its corresponding sequence of integrals $\left(
\int_{\left[  0,\text{ }\infty\right)  }\Gamma_{n}d\mu\right)  $\ been bounded
above by $C_{b_{1}}+C_{b_{2}}<\infty$, see Remark \ref{rem3}. Then by the
Beppo Levi's Theorem we can derive that sequence $\left(  \Gamma_{n}\right)  $
converges almost everywhere to some integrable measurable function $\Gamma$
and $\int_{\left[  0,\text{ }\infty\right)  }\Gamma d\mu=\lim_{n\rightarrow
\infty}\int_{\left[  0,\text{ }\infty\right)  }\Gamma_{n}d\mu\leq
\lim_{n\rightarrow\infty}\operatorname{d}\left(  \Delta_{n},\Phi_{2}\right)
=0$, meaning that $\lim_{n\rightarrow\infty}\inf_{k\geq n}\Delta_{k}=\Phi_{2}$
almost everywhere. There are two cases to be clarified. First assume that
$b_{1}<b_{2}$. Obviously, $\mu\left(  \left(  b_{1},~b_{2}\right)  \right)
>0$, so that there must be at least one point $x_{0}\in\left(  b_{1}%
,~b_{2}\right)  $ such that $\lim_{n\rightarrow\infty}\inf_{k\geq n}\Delta
_{k}\left(  x_{0}\right)  =\Phi_{2}\left(  x_{0}\right)  $. But since
$b_{1}<b_{2}$ the concave property implies that the graph of $\Phi_{2}%
$\ (resp. the graph of each function $\inf_{k\geq n}\Delta_{k}$) lies above
(resp. below) the graph of the line of equation $y=x$ in the interval $\left(
b_{1},~b_{2}\right)  $. Consequently, $\lim_{n\rightarrow\infty}\inf_{k\geq
n}\Delta_{k}\left(  x_{0}\right)  \leq x_{0}<\Phi_{2}\left(  x_{0}\right)  $.
This, however, is absurd since $\lim_{n\rightarrow\infty}\inf_{k\geq n}%
\Delta_{k}\left(  x_{0}\right)  =\Phi_{2}\left(  x_{0}\right)  $. Considering
the second case when $b_{1}>b_{2}$ we can similarly get into a contradiction.
Therefore, the statement is valid.
\end{proof}

\begin{lemma}
\label{lem4}Let $b$ and $c\in\left(  0,\,\infty\right)  $ be two arbitrary
numbers. Then the following assertions are equivalent:\newline\emph{(1)} The
equality $b=c$ holds.\newline\emph{(2)} The sets $\mathcal{Y}_{b}$ and
$\mathcal{Y}_{c}$ are equal.\newline\emph{(3)} The equality
$\operatorname*{dist}\left(  \mathcal{Y}_{b},\mathcal{Y}_{c}\right)  =0$ holds.
\end{lemma}

\begin{proof}
We first note that the chain of implications (1) $\rightarrow$ (2)
$\rightarrow$ (3) is obviously true. Thus we need only show the conditional
(3) $\rightarrow$ (1). In fact, assume that $\operatorname*{dist}\left(
\mathcal{Y}_{b},\mathcal{Y}_{c}\right)  =0$ but $b\neq c$. Then $\rho\left(
\mathcal{Y}_{b},\Delta\right)  =0$ for all $\Delta\in\mathcal{Y}_{c}$.
Nevertheless, this contradicts Lemma \ref{lem3}, since $b\neq c$. Therefore,
the argument is valid.
\end{proof}

We can similarly prove that:

\begin{lemma}
\label{lem5}Let $b$ and $c\in\left(  0,\,\infty\right)  $ be two arbitrary
numbers. Then the following assertions are equivalent\emph{:}\newline%
\emph{(1)} The equality $b=c$ holds.\newline\emph{(2)} The sets $\mathcal{A}%
_{b}$ and $\mathcal{A}_{c}$ are equal.\newline\emph{(3)} The equality
$\operatorname*{dist}\left(  \mathcal{A}_{b},\mathcal{A}_{c}\right)  =0$ holds.
\end{lemma}

\begin{theorem}
\label{theo6}Let $b$ and $c\in\left(  0,\,\infty\right)  $ be two arbitrary
numbers. Then the quantities $\operatorname*{dist}\left(  \mathcal{A}%
_{b},\mathcal{A}_{c}\right)  $\ and $\operatorname*{dist}\left(
\mathcal{Y}_{b},\mathcal{Y}_{c}\right)  $ define metrics on $\mathcal{Z}$ and
$\mathcal{Z}^{\ast}$ respectively. Hence, the couples $\left(  \mathcal{Z}%
,\operatorname*{dist}\right)  $ and $\left(  \mathcal{Z}^{\ast}%
,\operatorname*{dist}\right)  $ are metric spaces.
\end{theorem}

\begin{proof}
We need only show that $\operatorname*{dist}\left(  \mathcal{Y}_{b}%
,\mathcal{Y}_{c}\right)  $ is a metric on the set $\mathcal{Z}^{\ast}$,
because the other case can be similarly proved. In fact, we first point out
that the condition $\operatorname*{dist}\left(  \mathcal{Y}_{b},\mathcal{Y}%
_{c}\right)  \geq0$ is obvious and, by Lemma \ref{lem4} the equality holds if
and only if $\mathcal{Y}_{b}=\mathcal{Y}_{c}$.\ We also note that the symmetry
property trivially holds true. We are now left with the proof of the triangle
inequality. In fact, let $\mathcal{Y}_{b_{j}}\in\mathcal{Z}^{\ast}$ and
$\Phi_{j}\in\mathcal{Y}_{b_{j}}$\ ($j\in\left\{  1,~2,~3\right\}  $) be
arbitrary. Then by Proposition 5 (cf. \cite{AGB2005}, page 15) we have that
$\operatorname{d}\left(  \Phi_{1},\Phi_{3}\right)  \leq\operatorname{d}\left(
\Phi_{1},\Phi_{2}\right)  +\operatorname{d}\left(  \Phi_{2},\Phi_{3}\right)
$. Next, by taking the infimum over $\Phi_{3}\in\mathcal{Y}_{b_{3}}$ it
follows that
\[
\rho\left(  \Phi_{1},\mathcal{Y}_{b_{3}}\right)  \leq\operatorname{d}\left(
\Phi_{1},\Phi_{2}\right)  +\rho\left(  \Phi_{2},\mathcal{Y}_{b_{3}}\right)
\leq\operatorname{d}\left(  \Phi_{1},\Phi_{2}\right)  +\operatorname*{dist}%
\left(  \mathcal{Y}_{b_{2}},\mathcal{Y}_{b_{3}}\right)  ,
\]
i.e.\ $\rho\left(  \Phi_{1},\mathcal{Y}_{b_{3}}\right)  \leq\operatorname{d}%
\left(  \Phi_{1},\Phi_{2}\right)  +\operatorname*{dist}\left(  \mathcal{Y}%
_{b_{2}},\mathcal{Y}_{b_{3}}\right)  $. Finally, taking the infimum over
$\Phi_{2}\in\mathcal{Y}_{b_{2}}$\ yields $\rho\left(  \Phi_{1},\mathcal{Y}%
_{b_{3}}\right)  \leq\rho\left(  \Phi_{1},\mathcal{Y}_{b_{2}}\right)
+\operatorname*{dist}\left(  \mathcal{Y}_{b_{2}},\mathcal{Y}_{b_{3}}\right)
$, so that
\[
\operatorname*{dist}\left(  \mathcal{Y}_{b_{1}},\mathcal{Y}_{b_{3}}\right)
\leq\operatorname*{dist}\left(  \mathcal{Y}_{b_{1}},\mathcal{Y}_{b_{2}%
}\right)  +\operatorname*{dist}\left(  \mathcal{Y}_{b_{2}},\mathcal{Y}_{b_{3}%
}\right)  .
\]
This was to be proven.
\end{proof}

By the law of trichotomy it is not hard to see that $\left(  \mathcal{Z}%
,\preceq\right)  $ and $\left(  \mathcal{Z}^{\ast},\preceq\right)  $ are
lattices. Here too, the supremum and infimum binary operations on the lattices
$\left(  \mathcal{Z},\preceq\right)  $ and $\left(  \mathcal{Z}^{\ast}%
,\preceq\right)  $ will be denoted by the usual symbols $\vee$ and $\wedge$
respectively. We also point out that $\left(  \mathcal{Z},\preceq\right)  $
and $\left(  \mathcal{Z}^{\ast},\preceq\right)  $ are infinite graphs. Between
two vertices $\mathcal{A}_{b_{1}}$, $\mathcal{A}_{b_{2}}\in\mathcal{Z}$ we can
define the edge in two different ways: one by $e=\operatorname*{dist}\left(
\mathcal{A}_{b_{1}},\mathcal{A}_{b_{2}}\right)  \in\left(  0,~\infty\right)  $
and the other one by $\mathcal{A}_{e}\in\mathcal{Z}$ where
$e=\operatorname*{dist}\left(  \mathcal{A}_{b_{1}},\mathcal{A}_{b_{2}}\right)
$. These two edges can apply for the vertices of $\mathcal{Z}^{\ast}$ as well.

\section{Dense subsets in $\mathcal{Y}_{b}^{\left(  n\right)  }$}

\begin{theorem}
\label{theo7}Let $b\in\left(  0,\,\infty\right)  $ be an arbitrary number.
Then $\mathcal{A}_{b}$ is a dense set in $\mathcal{Y}_{b}$.
\end{theorem}

\begin{proof}
Fix arbitrarily any function $\Psi\in\mathcal{Y}_{b}$. Then there is some
$\Phi\in\mathcal{Y}_{\mathrm{conc}}$ such that $\Psi=\frac{b\Phi}{\Phi\left(
b\right)  }$ (by Lemma \ref{lem1}). Define $\Psi_{n}\left(  x\right)
=\frac{b\left(  \Phi\left(  x\right)  \right)  ^{1-1/(n+1)}}{\left(
\Phi\left(  b\right)  \right)  ^{1-1/(n+1)}}$, for all $x\in\left[
0,\,\infty\right)  $ and $n\in\mathbb{N}$. As we know from Theorem \ref{theo2}
(cf. \cite{AGB2005}, page 6) function $\Phi^{1-1/(n+1)}\in\mathcal{A}$ for all
$\Phi\in\mathcal{Y}_{\mathrm{conc}}$, $n\in\mathbb{N}$. Then $\left(  \Psi
_{n}\right)  \subset\mathcal{A}$ (via Lemma \ref{lem1}, \cite{AGB2005}, page
5). Hence, $\left(  \Psi_{n}\right)  \subset\mathcal{A}_{b}$, since $\Psi
_{n}\left(  b\right)  =b$ for all $n\in\mathbb{N}$. We can easily show that
$\left(  \Psi_{n}\right)  $ converges pointwise to $\Psi$. By Remark
\ref{rem4} it ensues that $\Psi\left(  x\right)  \leq h_{b}\left(  x\right)  $
and $\Psi_{n}\left(  x\right)  \leq h_{b}\left(  x\right)  $ for all
$x\in\left[  0,\,\infty\right)  $ and $n\in\mathbb{N}$, where $h_{b}\left(
x\right)  =x+b$, $x\in\left[  0,\,\infty\right)  $. We know via Remark
\ref{rem3} that function $h_{b}$ is square integrable. Then by applying twice
the Dominated Convergence Theorem one can verify that%
\[
\lim_{n\rightarrow\infty}%
{\displaystyle\int_{\left[  0,\,\infty\right)  }}
\Psi_{n}^{2}d\mu=%
{\displaystyle\int_{\left[  0,\,\infty\right)  }}
\Psi^{2}d\mu\text{ \ and }\lim_{n\rightarrow\infty}%
{\displaystyle\int_{\left[  0,\,\infty\right)  }}
\Psi_{n}\Psi d\mu=%
{\displaystyle\int_{\left[  0,\,\infty\right)  }}
\Psi^{2}d\mu,
\]
so that $\lim_{n\rightarrow\infty}\operatorname{d}\left(  \Psi,\Psi
_{n}\right)  =0$, because $\Psi\left(  x\right)  \Psi_{n}\left(  x\right)
\leq\left(  h_{b}\left(  x\right)  \right)  ^{2}$ for all $x\in\left[
0,\,\infty\right)  $ and $n\in\mathbb{N}$ (by Remark \ref{rem4}). This was to
be proven.
\end{proof}

\begin{theorem}
\label{theo8}Fix any pair of numbers $n\in\mathbb{N}\backslash\left\{
1\right\}  $ and $b\in\left(  0,~\infty\right)  $. Then $\mathcal{A}%
_{b}^{\left(  n\right)  }$ is dense in $\mathcal{Y}_{b}^{\left(  n\right)  }$.
\end{theorem}

\begin{proof}
Pick arbitrarily some $\Delta\in\mathcal{Y}_{b}^{\left(  n\right)  }$. Since
obviously $\mathcal{Y}_{b}^{CO\left(  n\right)  }$ is a proper subset of
$\mathcal{Y}_{b}^{\left(  n\right)  }$, we will have two cases to take into
consideration. First assume that $\Delta\in\mathcal{Y}_{b}^{CO\left(
n\right)  }$. This means that there can be found a counting number $k\geq n$
and a finite sequence $\Phi_{1},~\ldots,~\Phi_{k}\in\mathcal{Y}_{b}^{\left(
1\right)  }=\mathcal{Y}_{b}$ such that $\Delta=\Phi_{1}\circ\ldots\circ
\Phi_{k}$. Fix any integer $j\in\mathbb{N}$ and write $\Delta_{j}=\Psi
_{j}\circ\Delta$, where $\Psi_{j}\left(  x\right)  =\left(  b^{1/j}x\right)
^{j/(j+1)}$, $x\in\left[  0,~\infty\right)  $. Clearly, $\Psi_{j}%
\in\mathcal{A}_{b}^{\left(  1\right)  }$ for all $j\in\mathbb{N}$. Then
applying Theorem \ref{theo2} in \cite{AGB2005} and via the structure of set
$\mathcal{A}_{b}^{CO\left(  n\right)  }$, we can deduce that $\Delta_{j}%
\in\mathcal{A}_{b}^{CO\left(  n\right)  }$ for all $j\in\mathbb{N}$. It is not
difficult to see that sequence $\left(  \Delta_{j}\right)  $\ converge
pointwise to $\Delta$. By Remark \ref{rem4} we observe that $\Delta\leq h_{b}%
$, $\Delta_{j}\leq h_{b}$ and hence, $\Delta\Delta_{j}\leq\left(
h_{b}\right)  ^{2}$ on $\left[  0,~\infty\right)  $. Then recalling twice the
Dominated Convergence Theorem we can easily verify that
\[
\lim_{j\rightarrow\infty}\int_{\left[  0,~\infty\right)  }\left(  \Delta
_{j}\right)  ^{2}d\mu=\int_{\left[  0,~\infty\right)  }\Delta^{2}d\mu
=\lim_{j\rightarrow\infty}\int_{\left[  0,~\infty\right)  }\Delta\Delta
_{j}d\mu.
\]
Consequently, $\lim_{j\rightarrow\infty}\operatorname{d}\left(  \Delta
,\Delta_{j}\right)  =0$. In the second case we can suppose that $\Delta
\in\mathcal{Y}_{b}^{\left(  n\right)  }\backslash\mathcal{Y}_{b}^{CO\left(
n\right)  }$. Then without loss of generality we may choose $\Phi_{1}%
,~\ldots,~\Phi_{k}\in\mathcal{Y}_{b}^{CO\left(  n\right)  }$, whose graphs are
pairwise distinct, and some finite sequence $\left(  t_{1},~\ldots
~t_{k}\right)  \in H_{\left[  0,~1\right]  }$ with $\left(  t_{1}%
,~\ldots~t_{k}\right)  \subset\left(  0,~1\right)  $ such that $\Delta
=\sum_{i=1}^{k}t_{i}\Phi_{i}$. Consider $\Delta_{j}=\sum_{i=1}^{k}t_{i}\left(
\Psi_{j}\circ\Phi_{i}\right)  $, where $\Psi_{j}\left(  x\right)  =\left(
b^{1/j}x\right)  ^{j/(j+1)}$, $x\in\left[  0,~\infty\right)  $, $j\in
\mathbb{N}$. Clearly, on the one hand we have that $\left(  \Delta_{j}\right)
\subset\mathcal{A}_{b}^{\left(  n\right)  }$ because $\left(  \Psi_{j}%
\circ\Phi_{i}\right)  \subset\mathcal{A}_{b}^{CO\left(  n\right)  }$ for every
fixed index $i\in\left\{  1,~\ldots~k\right\}  $ and on the other hand
$\lim_{j\rightarrow\infty}\operatorname{d}\left(  \Phi_{i},\Psi_{j}\circ
\Phi_{i}\right)  =0$, $i\in\left\{  1,~\ldots~k\right\}  $, because of the
first part of this proof. Consequently, by the Minkowski inequality we can
observe that $\lim_{j\rightarrow\infty}\operatorname{d}\left(  \Delta
,\Delta_{j}\right)  \leq\sum_{i=1}^{k}t_{i}\lim_{j\rightarrow\infty
}\operatorname{d}\left(  \Phi_{i},\Psi_{j}\circ\Phi_{i}\right)  =0$. This
completes the proof.
\end{proof}

\section{Some criterium on the $L^{p}$-norm}

The result here below is worth being mentioned, which is an answer to the
second open problem in \cite{AGB2005}.

\begin{theorem}
\label{theo9}Let $\Phi\in\mathcal{Y}_{\mathrm{conc}}$ be arbitrary. Then the
following assertions are equivalent.\newline\emph{(1)} $\lim_{t\rightarrow
\infty}\frac{\Phi\left(  t\right)  }{t}=\lim_{t\rightarrow\infty}%
\varphi\left(  t\right)  \in\left(  0,\,\infty\right)  $.\newline\emph{(2)}
There is some constant $c\in\left[  1,\,\infty\right)  $ such that $c\Phi
>\Phi_{\operatorname{id}}$ on $\left(  0,\,\infty\right)  $.\newline\emph{(3)}
There is some constant $c\in\left[  1,\,\infty\right)  $ and some strictly
concave function $\Delta:\left[  0,\,\infty\right)  \rightarrow\left[
0,\,\infty\right)  $, differentiable on $\left(  0,\,\infty\right)  $ and
vanishing at the origin such that $c\Phi=\Phi_{\operatorname{id}}+\Delta$ on
$\left[  0,\,\infty\right)  $.
\end{theorem}

\begin{proof}
We first prove the conditional (1)$\rightarrow$(2). In fact, assume that
$\lim_{t\rightarrow\infty}\frac{\Phi\left(  t\right)  }{t}\in\left(
0,\,\infty\right)  $ but in the contrary for every counting number
$k\in\mathbb{N}$ there is some $x_{k}\in\left(  0,\,\infty\right)  $ for which
$k\Phi\left(  x_{k}\right)  \leq x_{k}$. Obviously, $\limsup
\limits_{k\rightarrow\infty}\frac{\Phi\left(  x_{k}\right)  }{x_{k}}\leq
\lim_{k\rightarrow\infty}k^{-1}=0$ which is absurd since $\limsup
\limits_{k\rightarrow\infty}\frac{\Phi\left(  x_{k}\right)  }{x_{k}}\in\left(
0,\,\infty\right)  $ by the assumption. Next we show the implication
(2)$\rightarrow$(3). In fact, assume that there is some constant $c\in\left[
1,\,\infty\right)  $ such that $c\Phi>\Phi_{\operatorname{id}}$ on $\left(
0,\,\infty\right)  $ and write $\Delta:=c\Phi-\Phi_{\operatorname{id}}$.
Clearly, $\Delta:\left[  0,\,\infty\right)  \rightarrow\left[  0,\,\infty
\right)  $ is a function such that $\Delta\left(  0\right)  =0$ and $\Delta$
is positive on $\left(  0,\,\infty\right)  $. We also note that $\Delta$ is
differentiable on $\left(  0,\,\infty\right)  $. Writing $\delta$ for the
derivative of $\Delta$, we can observe that $\delta=c\varphi-1$ on $\left(
0,\,\infty\right)  $. To show that $\Delta$ is strictly concave it is enough
if we prove that
\[
\left(  y-x\right)  \delta\left(  y-0\right)  <\Delta\left(  y\right)
-\Delta\left(  x\right)  <\left(  y-x\right)  \delta\left(  x+0\right)
=\left(  y-x\right)  \delta\left(  x\right)
\]
for all $x$, $y\in\left(  0,\,\infty\right)  $ with $x<y$ (where,
$\delta\left(  t-0\right)  $ respectively is the left derivative and
$\delta\left(  t+0\right)  $ the right derivative of $\Delta$ at point $t$).
In fact, fix arbitrarily two numbers $x$, $y\in\left(  0,\,\infty\right)  $
such that $x<y$. But since $\Phi$ is strictly concave we have that%
\[
\left(  y-x\right)  \varphi\left(  y-0\right)  <\Phi\left(  y\right)
-\Phi\left(  x\right)  <\left(  y-x\right)  \varphi\left(  x+0\right)
=\left(  y-x\right)  \varphi\left(  x\right)
\]
which easily leads to%
\[
c\varphi\left(  y-0\right)  <\frac{c\Phi\left(  y\right)  -c\Phi\left(
x\right)  }{y-x}<c\varphi\left(  x+0\right)  =c\varphi\left(  x\right)  .
\]
Hence,%
\[
c\varphi\left(  y-0\right)  -1<\frac{c\Phi\left(  y\right)  -c\Phi\left(
x\right)  }{y-x}-1<c\varphi\left(  x\right)  -1,
\]
i.e.
\[
\left(  y-x\right)  \delta\left(  y-0\right)  <\Delta\left(  y\right)
-\Delta\left(  x\right)  <\left(  y-x\right)  \delta\left(  x\right)  .
\]
This ends the proof of the implication (2)$\rightarrow$(3). In the last step,
we just point out that the conditional (3)$\rightarrow$(1) is obvious.
Therefore, we can conclude on the validity of the argument.
\end{proof}

Denote $\widetilde{\mathcal{Y}_{\mathrm{conc}}}:=\left\{  \Phi\in
\mathcal{Y}_{\mathrm{conc}}:\lim_{t\rightarrow\infty}\frac{\Phi\left(
t\right)  }{t}>0\right\}  $. It is not difficult to check that $\widetilde
{\mathcal{Y}_{\mathrm{conc}}}=\left\{  \Delta\in\mathcal{Y}_{\mathrm{conc}%
}:c\Delta>\Phi_{\operatorname{id}}\text{ on }\left(  0,\,\infty\right)  \text{
for some }c\in\left[  1,\,\infty\right)  \right\}  $. Write $T_{\Delta
}=\left\{  c\in\left[  1,\,\infty\right)  :c\Delta>\Phi_{\operatorname{id}%
}\text{ on }\left(  0,\,\infty\right)  \right\}  $, $\Delta\in\widetilde
{\mathcal{Y}_{\mathrm{conc}}}$.

Some few words about set $\widetilde{\mathcal{Y}_{\mathrm{conc}}}$.

\begin{remark}
\label{rem9}Let $\alpha\in\left(  0,\,\infty\right)  $ be arbitrary. Then
$\alpha\Delta\in\widetilde{\mathcal{Y}_{\mathrm{conc}}}$ provided that
$\Delta\in\widetilde{\mathcal{Y}_{\mathrm{conc}}}$.
\end{remark}

\begin{proof}
Whenever $\Delta\in\widetilde{\mathcal{Y}_{\mathrm{conc}}}$ we can choose a
corresponding $c\in T_{\Delta}$ such that $c\Delta>\Phi_{\operatorname{id}}$
on $\left(  0,\,\infty\right)  $. Now choose a constant $t_{0}\in\left(
1,\,\infty\right)  $ such that $\alpha t_{0}\geq c$. Hence, $t_{0}\left(
\alpha\Delta\right)  \geq c\Delta>\Phi_{\operatorname{id}}$ on $\left(
0,\,\infty\right)  $, i.e. $\alpha\Delta\in\widetilde{\mathcal{Y}%
_{\mathrm{conc}}}$.
\end{proof}

\begin{remark}
\label{rem10}Every function $\Delta\in\widetilde{\mathcal{Y}_{\mathrm{conc}}}$
can be written as the sum of a finite number of elements of $\widetilde
{\mathcal{Y}_{\mathrm{conc}}}$. Conversely, the sum of a finite number of
elements of $\widetilde{\mathcal{Y}_{\mathrm{conc}}}$\ also belongs to
$\widetilde{\mathcal{Y}_{\mathrm{conc}}}$.
\end{remark}

Next, we show that the quantities $\left\Vert f\right\Vert _{L^{p}}$ and
$\sup_{\Phi\in\widetilde{\mathcal{Y}_{\mathrm{conc}}}}\left(  \Phi\left(
1\right)  \right)  ^{-1}\left\Vert \Phi\circ\left\vert f\right\vert
\right\Vert _{L^{p}}$\ are equivalent, in the sense that they are both either
finite or infinite at the same time. This provides a kind of criterium for a
measurable function to belong to $L^{p}$.

\begin{theorem}
\label{theo10}Let $f$ be any measurable function on an arbitrarily fixed
measure space $\left(  \Omega,\mathcal{F},\lambda\right)  $ and $p\in\left[
1,\infty\right)  $ be any number. Then%
\[
\left\Vert f\right\Vert _{L^{p}}\leq\sup_{\Phi\in\widetilde{\mathcal{Y}%
_{\mathrm{conc}}}}\left(  \Phi\left(  1\right)  \right)  ^{-1}\left\Vert
\Phi\circ\left\vert f\right\vert \right\Vert _{L^{p}}\leq\left\Vert
f\right\Vert _{L^{p}}+\lambda\left(  \Omega\right)  .
\]

\end{theorem}

\begin{proof}
Pick any function $\Phi\in\widetilde{\mathcal{Y}_{\mathrm{conc}}}$. Then%
\[
\int_{\Omega}\left(  \left(  \Phi\left(  1\right)  \right)  ^{-1}\Phi
\circ\left\vert f\right\vert \right)  ^{p}d\lambda\leq\int_{\Omega}\left(
\left\vert f\right\vert +1\right)  ^{p}d\lambda
\]
because $\Delta\leq\left(  \Phi_{\operatorname{id}}+1\right)  \Delta\left(
1\right)  $ for all $\Delta\in\mathcal{Y}_{\mathrm{conc}}$. Consequently, via
the Minkowski inequality, it follows that $\left(  \Phi\left(  1\right)
\right)  ^{-1}\left\Vert \Phi\circ\left\vert f\right\vert \right\Vert _{L^{p}%
}\leq\left\Vert f\right\Vert _{L^{p}}+\lambda\left(  \Omega\right)  $, which
proves the inequality on the right hand-side of the above chain. To show the
left side inequality fix any $\Delta\in\mathcal{Y}_{\mathrm{conc}}$ and write
$\Delta_{n}=n^{-1}\Delta$, $n\in\mathbb{N}$. Clearly, $\left(  \Delta
_{n}\right)  \subset\mathcal{Y}_{\mathrm{conc}}$. It is also evident that
$\Phi_{\operatorname{id}}+\Delta_{n}\in\widetilde{\mathcal{Y}_{\mathrm{conc}}%
}$, $n\in\mathbb{N}$. Then
\begin{align*}
\sup_{\Phi\in\widetilde{\mathcal{Y}_{\mathrm{conc}}}}\left(  \Phi\left(
1\right)  \right)  ^{-1}\left\Vert \Phi\circ\left\vert f\right\vert
\right\Vert _{L^{p}}  &  \geq\left(  1+n^{-1}\right)  ^{-1}\left\Vert \left(
\Phi_{\operatorname{id}}+\Delta_{n}\right)  \circ\left\vert f\right\vert
\right\Vert _{L^{p}}\\
&  =\left(  1+n^{-1}\right)  ^{-1}\left\Vert \left\vert f\right\vert
+n^{-1}\left\vert f\right\vert \right\Vert _{L^{p}}\geq\left(  1+n^{-1}%
\right)  ^{-1}\left\Vert f\right\Vert _{L^{p}}.
\end{align*}
Passing to the limit yields $\sup_{\Phi\in\widetilde{\mathcal{Y}%
_{\mathrm{conc}}}}\left(  \Phi\left(  1\right)  \right)  ^{-1}\left\Vert
\Phi\circ\left\vert f\right\vert \right\Vert _{L^{p}}\geq\left\Vert
f\right\Vert _{L^{p}}$. Therefore, we have obtained a valid argument.
\end{proof}

\begin{theorem}
\label{theo11}Let $\left(  \Omega,\mathcal{F},\lambda\right)  $ be any measure
space and on it let $f$ be any measurable function. Then%
\[
\lambda\left(  \left\vert f\right\vert \geq\varepsilon\right)  =\inf\left\{
\inf\left\{  \lambda\left(  \Delta\circ\left\vert f\right\vert \geq\varepsilon
c^{-1}\right)  :c\in T_{\Delta}\right\}  :\Delta\in\widetilde{\mathcal{Y}%
_{\mathrm{conc}}}\right\}
\]
for every number $\varepsilon\in\left[  0,\,\infty\right)  $.
\end{theorem}

\begin{proof}
Throughout the proof $\varepsilon\in\left[  0,\,\infty\right)  $
will be any fixed number. We first note that the assertion is
trivial when $\left( \left\vert f\right\vert =\infty\right)
\neq\varnothing$. We shall then prove it when $\left(  \left\vert
f\right\vert <\infty\right)  \neq\varnothing$. Pick some
$\Delta\in\widetilde{\mathcal{Y}_{\mathrm{conc}}}$ and $c\in
T_{\Delta}$ such that $c\Delta>\Phi_{\operatorname{id}}$ on
$\left( 0,\,\infty\right)  $. It is not hard to see that $\left(
\left\vert f\right\vert \geq\varepsilon\right)  =\left(
\Delta\circ\left\vert f\right\vert \geq\Delta\left(
\varepsilon\right)  \right)  \subset\left(
\Delta\circ\left\vert f\right\vert \geq\varepsilon c^{-1}\right)  $ and thus%
\[
\lambda\left(  \left\vert f\right\vert \geq\varepsilon\right)  =\lambda\left(
\Delta\circ\left\vert f\right\vert \geq\Delta\left(  \varepsilon\right)
\right)  \leq\lambda\left(  \Delta\circ\left\vert f\right\vert \geq\varepsilon
c^{-1}\right)  .
\]
Consequently,
\[
\lambda\left(  \left\vert f\right\vert \geq\varepsilon\right)  \leq
\inf\left\{  \inf\left\{  \lambda\left(  \Delta\circ\left\vert f\right\vert
\geq\varepsilon c^{-1}\right)  :c\in T_{\Delta}\right\}  :\Delta\in
\widetilde{\mathcal{Y}_{\mathrm{conc}}}\right\}  .
\]
To prove the converse statement, we need show that
\[
\lambda\left(  \left\vert f\right\vert \geq\varepsilon\right)  \geq
\inf\left\{  \inf\left\{  \lambda\left(  \Delta\circ\left\vert f\right\vert
\geq\varepsilon c^{-1}\right)  :c\in T_{\Delta}\right\}  :\Delta\in
\widetilde{\mathcal{Y}_{\mathrm{conc}}}\right\}  .
\]
In fact, for any $n\in\mathbb{N}$ set $\Delta_{n}=\Phi_{\operatorname{id}%
}+n^{-1}\left(  1-e^{-\Phi_{\operatorname{id}}}\right)  $. It is not difficult
to see that $\Delta_{n}\in\mathcal{Y}_{\mathrm{conc}}$ and $\Delta_{n}%
>\Phi_{\operatorname{id}}$ on $\left(  0,\,\infty\right)  $, $n\in\mathbb{N}$.
This means that $\left(  \Delta_{n}\right)  \subset\widetilde{\mathcal{Y}%
_{\mathrm{conc}}}$ and moreover, $1\in T_{\Delta_{n}}$, $n\in\mathbb{N}$.
Consequently,
\begin{align*}
\inf\left\{  \inf\left\{  \lambda\left(  \Delta\circ\left\vert f\right\vert
\geq\varepsilon c^{-1}\right)  :c\in T_{\Delta}\right\}  :\Delta\in
\widetilde{\mathcal{Y}_{\mathrm{conc}}}\right\}   &  \leq\lambda\left(
\Delta_{n}\circ\left\vert f\right\vert \geq\varepsilon\right)  \\
&  =\lambda\left(  \left\vert f\right\vert +n^{-1}\left(  1-e^{-\left\vert
f\right\vert }\right)  \geq\varepsilon\right)  \text{.}%
\end{align*}
However, as $\left(  \Delta_{n}\right)  $ is a decreasing sequence it is
obvious that $\left(  \Delta_{n+1}\circ\left\vert f\right\vert \geq
\varepsilon\right)  \subset\left(  \Delta_{n}\circ\left\vert f\right\vert
\geq\varepsilon\right)  $, $n\in\mathbb{N}$. Thus having passed to the limit
we can observe that
\[
\inf\left\{  \inf\left\{  \lambda\left(  \Delta\circ\left\vert f\right\vert
\geq\varepsilon c^{-1}\right)  :c\in T_{\Delta}\right\}  :\Delta\in
\widetilde{\mathcal{Y}_{\mathrm{conc}}}\right\}  \leq\lambda\left(  \left\vert
f\right\vert \geq\varepsilon\right)  .
\]
Therefore, the proof is a valid argument.
\end{proof}

\begin{theorem}
\label{theo12}Let $f\in L^{p}\left(  \Omega,\mathcal{F},\lambda\right)  $,
$p\geq1$, where $\left(  \Omega,\mathcal{F},\lambda\right)  $ is any given
measure space. Then
\[
\left\Vert f\right\Vert _{L^{p}}=\inf\left\{  \inf\left\{  c\left\Vert
\Delta\circ\left\vert f\right\vert \right\Vert _{L^{p}}:c\in T_{\Delta
}\right\}  :\Delta\in\widetilde{\mathcal{Y}_{\mathrm{conc}}}\right\}  .
\]

\end{theorem}

\begin{proof}
Pick arbitrarily some $\Delta\in\widetilde{\mathcal{Y}_{\mathrm{conc}}}$ and
$c\in T_{\Delta}$ such that $c\Delta>\Phi_{\operatorname{id}}$ on $\left(
0,\,\infty\right)  $. Clearly, $c\left\Vert \Delta\circ\left\vert f\right\vert
\right\Vert _{L^{p}}\geq\left\Vert f\right\Vert _{L^{p}}$. We can then easily
observe that
\[
\inf\left\{  \inf\left\{  c\left\Vert \Delta\circ\left\vert f\right\vert
\right\Vert _{L^{p}}:c\in T_{\Delta}\right\}  :\Delta\in\widetilde
{\mathcal{Y}_{\mathrm{conc}}}\right\}  \geq\left\Vert f\right\Vert _{L^{p}}.
\]
To prove the converse of this inequality consider the sequence $\left(
\Delta_{n}\right)  \subset\widetilde{\mathcal{Y}_{\mathrm{conc}}}$, where
$\Delta_{n}=\Phi_{\operatorname{id}}+n^{-1}\left(  1-e^{-\Phi
_{\operatorname{id}}}\right)  >\Phi_{\operatorname{id}}$ on $\left(
0,\,\infty\right)  $, $n\in\mathbb{N}$. Then as $1\in T_{\Delta_{n}}$,
$n\in\mathbb{N}$, we have%
\[
\inf\left\{  \inf\left\{  c\left\Vert \Delta\circ\left\vert f\right\vert
\right\Vert _{L^{p}}:c\in T_{\Delta}\right\}  :\Delta\in\widetilde
{\mathcal{Y}_{\mathrm{conc}}}\right\}  \leq\left\Vert \Delta_{n}%
\circ\left\vert f\right\vert \right\Vert _{L^{p}}.
\]
Since $\left(  \Delta_{n}\right)  $ is a decreasing sequence it ensues
that\ $\left(  \Delta_{n}\circ\left\vert f\right\vert \right)  $ is also a
decreasing sequence which tends to $\left\vert f\right\vert $. As every member
of sequence $\left(  \Delta_{n}\circ\left\vert f\right\vert \right)  $ is
dominated by $\Delta_{1}\circ\left\vert f\right\vert \in L^{p}$, then by
applying the Dominated Convergence Theorem it will entail that
\[
\inf\left\{  \inf\left\{  c\left\Vert \Delta\circ\left\vert f\right\vert
\right\Vert _{L^{p}}:c\in T_{\Delta}\right\}  :\Delta\in\widetilde
{\mathcal{Y}_{\mathrm{conc}}}\right\}  \leq\left\Vert f\right\Vert _{L^{p}}.
\]
This completes the proof.
\end{proof}

\begin{corollary}
\label{cor3}Suppose that $h:\mathbb{R\rightarrow R}$ is a continuous function.
Then $\left\vert h\right\vert =\frac{1}{\sqrt[p]{\lambda\left(  \Omega\right)
}}\inf\left\{  \inf\left\{  \left(  \Delta\circ\left\vert h\right\vert
\right)  c:c\in T_{\Delta}\right\}  :\Delta\in\widetilde{\mathcal{Y}%
_{\mathrm{conc}}}\right\}  $.
\end{corollary}

\begin{proof}
Fix any number $x\in\mathbb{R}$ and let $f\in L^{p}\left(  \Omega
,\mathcal{F},\lambda\right)  $\ be the constant function defined by $f\equiv
h\left(  x\right)  $\ on $\Omega$. Then by applying Theorem \ref{theo10} we
can easily deduce the result.
\end{proof}

\begin{openproblem}
Given any number $k\in\mathbb{N}$ characterize all pairs of functions $\Phi$
and $\Delta\in\mathcal{Y}_{\mathrm{conc}}$ such that $\left\vert \left\{
x\in\left(  0,~\infty\right)  :\Phi\left(  x\right)  =\Delta\left(  x\right)
\right\}  \right\vert =k$.
\end{openproblem}

\begin{openproblem}
Characterize all pairs of functions $\Phi$ and $\Delta\in\mathcal{Y}%
_{\mathrm{conc}}$ such that the sets $\left(  0,~\infty\right)  $ and
$\left\{  x\in\left(  0,~\infty\right)  :\Phi\left(  x\right)  =\Delta\left(
x\right)  \right\}  $ should be equinumerous.
\end{openproblem}

\end{document}